\documentclass[article, 1p, sort&compress]{elsarticle}

\usepackage{lineno,hyperref}
\usepackage{amssymb}
\usepackage{amsmath}
\usepackage{amsthm}

\theoremstyle{remark}
\newtheorem{remark}{Remark}

\theoremstyle{plain}
\newtheorem{proposition}{Proposition}

\theoremstyle{plain}
\newtheorem{lemma}{Lemma}

\theoremstyle{plain}
\newtheorem{theorem}{Theorem}

\theoremstyle{plain}

\theoremstyle{plain}
\newtheorem{corollary}{Corollary}

\newcommand{\Real}{\mathbb R}

\newcommand{\real}[1]{{\mathbb R}^{#1}}

\newcommand{\bu}{{\boldsymbol u}}

\newcommand{\bv}{{\boldsymbol v}}

\newcommand{\bx}{\boldsymbol x}
\newcommand{\by}{{\boldsymbol y}}

\newcommand{\buf}{{\bu(\cdot)}}  

\newcommand{\bzero}{{\bf 0}}

\newcommand{\bbeta}{{\mbox{\boldmath $\beta$}}}

\newcommand{\blam}{{\mbox{\boldmath $\lambda$}}}

\modulolinenumbers[65]

\journal{the Journal of \LaTeX\ Templates}

\begin{document}
\nolinenumbers
\begin{frontmatter}

\title{A Derivation of Nesterov's Accelerated Gradient Algorithm from Optimal Control Theory}


\author{I. M. Ross\fnref{myfootnote}}
\address{Naval Postgraduate School, Monterey, CA 93943}
\fntext[myfootnote]{Distinguished Professor \& Program Director, Control and Optimization }


\begin{abstract}
Nesterov's accelerated gradient algorithm is derived from first principles. The first principles are founded on the recently-developed optimal control theory for optimization. This theory frames an optimization problem as an optimal control problem whose trajectories generate various continuous-time algorithms.  The algorithmic trajectories satisfy the necessary conditions for optimal control.  The necessary conditions produce a controllable dynamical system for accelerated optimization.
Stabilizing this system via a quadratic control Lyapunov function generates an ordinary differential equation. An Euler discretization of the resulting differential equation produces Nesterov's algorithm. In this context, this result solves the purported mystery surrounding the algorithm.
 \\
\end{abstract}


\end{frontmatter}


\section{Introduction}
In broad terms, Nesterov's accelerated gradient method for minimizing a convex function, $E:\real{N_x} \to \Real$, is given by\cite{nesterov83},
\begin{subequations}\label{eq:NAG}
\begin{align}
\bx_k &= \by_k -\alpha_k\, \partial_{\by}E(\by_k) \label{eq:step-1:nes}\\
\by_{k+1} & = \bx_k + \beta_k (\bx_k - \bx_{k-1}) \label{eq:step-2:nes}
\end{align}
\end{subequations}
where, $\bx_k \in \real{N_x}, \ \by_k \in \real{N_x},\ \alpha_k \in \Real_+,\ \beta_k \in \Real_+$ and $ k \in \mathbb{N}$.  There are many ways to ``explain'' this algorithm starting from a well-deserved attribution to Nesterov's penetrating insights in convex programming to discretizations of certain ordinary differential equations(ODEs)\cite{shi:hi-res-Nes,alvarez2002,su}.  The use of ODEs to explain algorithms has a long and rich history\cite{boggs71,smale76,brown+biggs} going all the way back to Gavurin's pioneering work in the late 1950s\cite{gavurin}.  The ODEs that explain algorithms are typically derived by considering the limiting cases of the algorithmic maps themselves. In other words, an ODE is usually generated after an algorithm is invented and not the other way around.  Thus the question remains: can the ODEs be generated by some higher-level universal principle and without any a priori knowledge of algorithms?  The recently-developed optimal control theory for optimization\cite{rossJCAM-1} answers this question in the affirmative. According to this theory, controllable ODEs for optimization algorithms can be generated as outcomes of the necessary conditions for optimal control.  Producing a practical algorithm is then reduced to using suitable semi-discretization methods. As shown in \cite{rossJCAM-1}, these ideas generate a vast number of well-known but \emph{unaccelerated} algorithms such as Newton's method and the gradient method.  It was conjectured in \cite{rossJCAM-1} that the theory could also generate accelerated optimization methods by switching the dynamical model from a single integrator to a double integrator.  In this note, we prove this conjecture by deriving Nesterov's accelerated gradient method using optimal control theory as a foundation for optimization.\footnote{All the results presented in this note were initially contained in an earlier draft of \cite{ross-accelerated-arxiv}.  Because the current version of \cite{ross-accelerated-arxiv} contains a more comprehensive theory for accelerated optimization than its earlier counterparts, the ``Nesterov section'' was eliminated in the later drafts in favor of continuity of ideas.  In this context, this technical note is simply a parsed out version of an earlier draft of \cite{ross-accelerated-arxiv}. }

\section{Background: Optimal Control Theory for Optimization}
We rely heavily on \cite{rossJCAM-1} in developing the basic framework while providing sufficient details for completeness. To this end, consider the unconstrained static optimization
problem given by,
\begin{equation}
(S) \left\{\displaystyle\mathop\text{Minimize }_{\bx_f \in \real{N_x}}   E(\bx_f) \right.
\end{equation}
To produce an optimal control problem that solves $(S)$, we create a vector field by ``sweeping'' the function $E$ backwards in time according to,
\begin{equation}\label{eq:idea-1}
y(t) := E(\bx(t))
\end{equation}
Differentiating \eqref{eq:idea-1} with respect to time we get,
\begin{equation}\label{eq:cost-evolution}
\dot y = \partial_{\bx} E(\bx)\cdot \dot\bx
\end{equation}
As shown in \cite{rossJCAM-1}, if we set $\dot\bx = \bu$ to generate a controllable dynamical system, the resulting theory generates unaccelerated optimization methods.   In pursuit of acceleration, we replace the single integrator model by a double integrator,
\begin{equation}\label{eq:dint-statespace}
\dot\bx = \bv, \quad \dot\bv = \bu
\end{equation}
This model generates a primal controllable dynamical system given by,
\begin{equation}\label{eq:cost-evolution}
\dot y = \partial_{\bx} E(\bx)\cdot\bv, \quad \dot\bx = \bv, \quad \dot\bv = \bu
\end{equation}
Taking generic initial conditions and a final rest ``velocity,'' $\bv(t_f) = \bzero$, as boundary conditions, we arrive at the following candidate optimal control problem $(R)$ that purportedly solves the optimization problem $(S)$:
\begin{eqnarray}
&(R) \left\{
\begin{array}{lrl}
\textsf{Minimize }  & J[y(\cdot), \bx(\cdot), \bv(\cdot), \buf, t_f]
:=& y(t_f)  \\
\textsf{Subject to} & \dot\bx=& \bv\\
&\dot\bv = & \bu \\
&\dot y=& \partial_{\bx} E(\bx)\cdot \bv\\
&(\bx(t_0), t_0) =& (\bx^0, t^0) \\
& y(t_0) = & E(\bx^0) \\
&\bv(t_f) = & \bzero
\end{array} \right.& \label{eq:prob-R}
\end{eqnarray}
where, $\bx^0$ is an initial ``guess'' of the solution (to Problem $(S)$).  The variables $t_f, \bx(t_f)$ and $\bv(t_0)$ are all free.
\begin{remark}
The cost functional in Problem $(R)$ is given by the final value of the $y$ variable, which, by construction, is exactly equal to the objective function of Problem~$(S)$. Consequently, a solution to Problem $(R)$ generates a solution to Problem $(S)$.
\end{remark}
\begin{remark}
A solution to Problem $(R)$ generates an optimal $\bx$-trajectory.  A discretization of this continuous-time trajectory generates a practical algorithm for Problem $(S)$.
\end{remark}
It follows from the preceding remarks that not only is Problem $(S)$ embedded in Problem $(R)$  but also that a solution to Problem $(R)$ automatically generates a continuous-time algorithm for solving Problem $(S)$.

\section{Necessary Conditions for Problem $(R)$}

%
\begin{lemma}\label{lemma:normality}
Problem $(R)$ has no abnormal extremals.
\end{lemma}
\begin{proof}
The Pontryagin Hamiltonian\cite{ross-book} for this problem is given by,
\begin{equation}\label{eq:H-unc}
H(\blam_x, \blam_v, \lambda_y, \bx, \bv, y, \bu):= \blam_x \cdot \bv + \blam_v \cdot \bu + \lambda_y\, \partial_{\bx} E(\bx)\cdot \bv
\end{equation}
where, $\blam_{x}, \blam_v$ and $\lambda_y$ are costates that satisfy the adjoint equations,
\begin{subequations}
\begin{align}
\dot\blam_{x} &=-\partial_{\bx}H = -\lambda_y\,\partial^2_{\bx}E(\bx)\, \bv \label{eq:adj-x-unc}\\
\dot\blam_v & = -\partial_{\bv}H = - \blam_x -\lambda_y\,\partial_{\bx}E(\bx)\label{eq:adj-v-unc}\\
\dot\lambda_y & =-\partial_y H = 0 \label{eq:adj-y-unc}
\end{align}
\end{subequations}
The transversality conditions\cite{ross-book} for Problem $(R)$ are given by,
\begin{subequations}
\begin{align}
\blam_x(t_f) &= \bzero \label{eq:tvc-x}\\
\blam_v(t_0) & = \bzero \label{eq:tvc-v}\\
\blam_y(t_f) & = \nu_0 \ge 0 \label{eq:tvc-y}
\end{align}
\end{subequations}
where, $\nu_0$ is the cost multiplier. From \eqref{eq:adj-y-unc} and \eqref{eq:tvc-y} we have,
\begin{equation}\label{eq:lam-y-pre-tvc}
\lambda_y(t) = \nu_0
\end{equation}
If $\nu_0 = 0$, then $\lambda_y(t) \equiv 0$.  This implies, from \eqref{eq:adj-x-unc} and \eqref{eq:tvc-x}, that $\blam_x(t) \equiv \bzero$.  Similarly, $\blam_v(t) \equiv \bzero$ from \eqref{eq:adj-v-unc} and \eqref{eq:tvc-v}.  The vanishing of all multipliers violates the nontriviality condition.  Hence $\nu_0 > 0$.
\end{proof}
%
%
%
\begin{theorem}\label{thm-singular}
All extremals of Problem $(R)$ are singular.  Furthermore, the singular arcs are of infinite order.
\end{theorem}
\begin{proof}
The Hamiltonian is linear in the control variable and the control space is unbounded; hence, if $\bu$ is optimal, it must be singular.  Furthermore, from the Hamiltonian minimization condition we have the first-order condition,
\begin{equation}\label{eq:hmc-order-1}
\partial_{\bu} H = \blam_v(t) = \bzero \qquad\forall t \in [t_0, t_f]
\end{equation}
Differentiating \eqref{eq:hmc-order-1} with respect to time, we get,
\begin{equation}\label{eq:SA-deg-1}
\frac{d}{dt}\partial_{\bu} H =\dot\blam_v(t) = - \blam_x -\nu_0\,\partial_{\bx}E(\bx)  = \bzero
\end{equation}
Equation \eqref{eq:SA-deg-1} does not generate an expression for the control function; hence, taking the second time derivative of $\partial_{\bu}H$ we get,
\begin{align}
\frac{d^2}{dt^2}\partial_{\bu} H &= -\dot\blam_x - \nu_0\, \partial^2_{\bx} E(\bx)\,\dot\bx   \nonumber\\
        &=-\dot\blam_x - \nu_0\, \partial^2_{\bx} E(\bx)\,\bv \nonumber \\
        &\equiv \bzero
\end{align}
where, the last equality follows from \eqref{eq:adj-x-unc} and Lemma \ref{lemma:normality}.  Hence, we have,
$$ \frac{d^k}{dt^k}\partial_{\bu} H = \bzero \quad\text{for\ } k = 0, 1 \ldots  $$
and no $k$ yields an expression for $\bu$.
\end{proof}
%

\begin{theorem}[A Transversality Mapping Theorem]\label{thm:tmt}
The first-order necessary condition for Problem~$(S)$ is embedded in the terminal transversality condition for Problem $(R)$.
\end{theorem}
\begin{proof}
From \eqref{eq:SA-deg-1}, we have
\begin{equation}\label{eq:lam-x-sol-generic}
\blam_x(t) = -\nu_0\, \partial_{\bx} E(\bx(t))
\end{equation}
From \eqref{eq:tvc-x} and Lemma 1, it follows that $\partial_{\bx}E(\bx_f) = \bzero$.
\end{proof}
Collecting all relevant equations, it follows that the primal-dual control dynamical system generated by Problem $(R)$ is given by,
\begin{subequations}\label{eq:p-d-dynamics}
\begin{align}
\dot\bx & = \bv         & \dot\blam_{x} & = -\lambda_y\,\partial^2_{\bx}E(\bx) \, \bv \\
\dot\bv & = \bu         & \dot\blam_v & =  -\blam_x - \lambda_y\, \partial_{\bx} E(\bx)\\
\dot y & = \partial_{\bx} E(\bx)\cdot \bv    & \dot\lambda_y & = 0
\end{align}
\end{subequations}
The boundary conditions for \eqref{eq:p-d-dynamics} are given by,
\begin{subequations}\label{eq:BCs}
\begin{align}
\bx(t^0) & = \bx^0  & \bv(t_f) & = \bzero\\
y(t^0) &= E(\bx^0)    & \blam_x(t_f) & = \bzero\\
\blam_v(t^0) & =  \bzero  & \lambda_y(t_f) &= \nu_0 > 0
\end{align}
\end{subequations}
Because the optimal control is a singular arc of infinite order, an optimal trajectory of Problem~$(R)$ must satisfy \eqref{eq:p-d-dynamics} and \eqref{eq:BCs}. Along a singular arc, $\blam_v(t) \equiv \bzero$; hence, the auxiliary controllable dynamical system of interest\cite{rossJCAM-1} resulting from \eqref{eq:p-d-dynamics} is given by,
\begin{equation}\label{eq:aux-dynamics}
(A) \left \{
\begin{aligned}
\dot\blam_{x} & = -\partial^2_{\bx}E(\bx) \, \bv   \\
\dot\bv & = \bu
\end{aligned}
\right.
\end{equation}
where, we have scaled the adjoint covector $\blam_x$ by $\nu_0 > 0$ (cf. Lemma \ref{lemma:normality}).  The final-time condition for $(A)$ is given by,
\begin{equation}\label{eq:aux-target}
(T) \left \{
\begin{aligned}
\blam_x(t_f) & = \bzero\\
\bv(t_f) & = \bzero
\end{aligned}
\right.
\end{equation}
That is, any singular control that satisfies \eqref{eq:aux-dynamics} and \eqref{eq:aux-target} generates a candidate ``optimal'' continuous-time algorithm for Problem $(S)$.

\section{A Feedback Controller for the $(A)$-$(T)$ System}

Let $\bbeta$ be the control vector field defined according to,
\begin{equation}
\bbeta(\bx, \bv, \bu) := \left[
                        \begin{array}{c}
                          -\partial^2_{\bx}E(\bx) \, \bv\\
                           \bu\\
                        \end{array}
                      \right]
\end{equation}
Let $V:(\blam_x, \bv) \mapsto \Real$ be a control Lyapunov function (CLF)\cite{sontag-book} for the $(A, T)$ pair.
Let $\pounds_\beta V$ be the Lie derivative of $V$ along the vector field $\bbeta$.
Then, a sufficient condition\cite{sontag-book,clarkeLyap} for globally guiding the pair $(\blam_x, \bv)$ to $(\bzero, \bzero)$ is to render $\pounds_\beta V$ negative; i.e., we need to find a $\bu$ such that,
\begin{equation}\label{eq:clf-theory-1}
\pounds_\beta V := \partial V(\blam_x, \bv)\cdot \bbeta(\bx, \bv, \bu) < 0
\end{equation}
whenever $\bx \ne \bx_f$. Consequently, we seek to design a (singular) control function that satisfies \eqref{eq:clf-theory-1}.

For the remainder of this note, we choose the following positive definite CLF,
\begin{equation}\label{eq:V-quad-pd}
V(\blam_x, \bv) =  (a/2)\blam_x \cdot \blam_x + (b/2)\bv \cdot \bv+ c\blam_x \cdot \bv
\end{equation}
where,
\begin{equation}
 a > 0, \quad b > 0, \quad c < 0, \quad ab - c^2 > 0
\end{equation}
are constants. As a result, we have,
\begin{equation}\label{eq:Vdot-quadV}
\pounds_\beta V = -[a\blam_x + c\bv]\cdot \partial^2_{\bx}E(\bx)\, \bv + [c \blam_x + b \bv] \cdot \bu
\end{equation}
A generic linear feedback controller\cite{sontag-book} is given by $\bu = K_a\, \blam_x + K_b\, \bv$, where $K_a$ and $K_b$ are real numbers.  Motivated by the intuition to design a control that directly incorporates the drift vector field to render $\pounds_\beta V < 0$, consider a modification to the linear feedback control strategy given by,
\begin{equation}\label{eq:uFamilyNo2}
\bu = K_a\, \blam_x + K_b\, \bv + K_c \, \partial^2_{\bx}E(\bx)\, \bv
\end{equation}
where $K_a, K_b$ and $K_c$ are all real numbers that must be chosen to guarantee $\pounds_\beta V$ negative.
%
\begin{proposition}
Suppose $E$ is a convex function and $\bu$ is given by \eqref{eq:uFamilyNo2}.   If
\begin{equation}
K_a > 0, \quad K_b < 0, \quad b K_a  = c K_b,    \quad \text{and}\quad K_c = a/c
\end{equation}
then, $\pounds_\beta V < 0$ for all $(\blam_x, \bv) \ne \bzero$.
\end{proposition}
\begin{proof}
Substituting \eqref{eq:uFamilyNo2} in \eqref{eq:Vdot-quadV} we get,
\begin{eqnarray}\label{eq:Vdot-wu2}
\pounds_\beta V &= \big(-a + c K_c\big)\blam_x \cdot \partial^2_{\bx}E(\bx)\, \bv + \big(-c + b K_c\big) \bv \cdot \partial^2_{\bx}E(\bx)\, \bv \nonumber\\
 &  + \big(c \blam_x + b \bv\big)\cdot \big(K_a\, \blam_x + K_b\, \bv\big)
\end{eqnarray}
Substituting $b K_a = c K_b$ in the third term of \eqref{eq:Vdot-wu2} generates,
\begin{equation}\label{eq:term3-result}
\big(c \blam_x + b \bv\big) \cdot \big(K_a\, \blam_x + K_b\, \bv\big) =\frac{K_b}{b}\big(c \blam_x + b \bv\big) \cdot \big(c \blam_x + b \bv\big) \le 0
\end{equation}
where, the inequality in \eqref{eq:term3-result} follows from the assumption that $K_b < 0$.

With $K_c = a/c$, the first term of \eqref{eq:Vdot-wu2} vanishes.  The second term simplifies to,
\begin{equation}
\big(-c + b K_c\big) \bv \cdot \partial^2_{\bx}E(\bx)\, \bv = \left(\frac{-c^2 + ab}{c} \right)\bv \cdot \partial^2_{\bx}E(\bx)\, \bv
\end{equation}
Because $ab - c^2 >0$ and $c <0$, it follows that the second term of \eqref{eq:Vdot-wu2} is negative for a positive definite Hessian; hence, $\pounds_\beta V < 0$.
\end{proof}
%
\begin{corollary}
Let,
\begin{subequations}\label{eq:K2gamma}
\begin{align}
K_a &:= \gamma_a, &\gamma_a > 0\\
K_b &:= -\gamma_b, & \gamma_b > 0 \\
K_c &:= - \gamma_c, & \gamma_c > 0
\end{align}
\end{subequations}
then, the singular control law given by \eqref{eq:uFamilyNo2} generates the second order ODE,
\begin{equation}\label{eq:Nesterov}
\ddot\bx + \gamma_a\, \partial_{\bx} E(\bx) + \gamma_b\, \dot\bx + \gamma_c \, \partial^2_{\bx}E(\bx)\, \dot\bx = \bzero
\end{equation}
\end{corollary}
\begin{proof}
Equation \eqref{eq:Nesterov} follows directly from \eqref{eq:dint-statespace}, \eqref{eq:uFamilyNo2} and \eqref{eq:K2gamma}.
\end{proof}
%

\section{Equation \eqref{eq:Nesterov} Generates \eqref{eq:NAG}}

As shown by Shi et al\cite{shi:hi-res-Nes}, a discretization of \eqref{eq:Nesterov} generates Nesterov's accelerated gradient method.  To see this, consider first a discretization of the the last term on the left-hand-side of \eqref{eq:Nesterov}:
\begin{equation}\label{eq:term-3}
\gamma_c \, \partial^2_{\bx}E(\bx)\, \dot\bx = \gamma_c \,\frac{d}{dt}\Big(\partial_{\bx}E(\bx)\Big) \longrightarrow \frac{\gamma_c}{h_k}\Big(\partial_{\bx}E(\bx_k) - \partial_{\bx}E(\bx_{k-1})\Big)
\end{equation}
where, $h_k > 0$ is a discretization step. Next, consider the first three terms of \eqref{eq:Nesterov}. These are identical to Polyak's equation whose discretization generates the heavy ball method\cite{polyak64,polyak-ode},
\begin{equation}\label{eq:polyak}
\bx_{k+1} = \bx_k - \alpha_k \partial_{\bx}E(\bx_k) + \beta_k (\bx_k - \bx_{k-1})
\end{equation}
Hence, \eqref{eq:Nesterov} may be discretized as,
\begin{equation}\label{eq:Nesterov-derived}
\bx_{k+1} = \bx_k - \alpha_k \partial_{\bx}E(\bx_k) + \beta_k (\bx_k - \bx_{k-1}) - \gamma_k \Big(\partial_{\bx}E(\bx_k) - \partial_{\bx}E(\bx_{k-1})\Big)
\end{equation}
Substituting \eqref{eq:step-1:nes} in \eqref{eq:step-2:nes}, Nesterov's method for $\alpha_k = \alpha$ may be rewritten as,
\begin{equation}\label{eq:Nesterov-y}
\by_{k+1} = \by_k - \alpha\, \partial_{\by}E(\by_k) + \beta_k (\by_k - \by_{k-1}) - \alpha\,\beta_k \Big(\partial_{\by}E(\by_k) - \partial_{\by}E(\by_{k-1})\Big)
\end{equation}
Equation \eqref{eq:Nesterov-y} is the same as \eqref{eq:Nesterov-derived} with $\gamma_k = \alpha\beta_k$.

\begin{remark}
Equation \eqref{eq:Nesterov} was introduced and studied by Alvarez et al\cite{alvarez2002} as a ``dynamical inertial Newton'' system.  Shi et al \cite{shi:hi-res-Nes} generated this system as a ``high-resolution'' ODE that represents Nesterov's  method\cite{nesterov83}.  This ODE (i.e., \eqref{eq:Nesterov}) is different from the one generated in \cite{su} to model Nesterov's method. The ODE in \cite{su} does not contain the Hessian term (i.e., $\gamma_c = 0$) but has $\gamma_b$ time varying.  An ODE with time-varying gains (i.e., $\gamma_a(t), \gamma_b(t), \gamma_c(t)$) can be generated using the more general method developed in \cite{ross-accelerated-arxiv}.
\end{remark}
\begin{remark}
The condition $\alpha_k = \alpha$ used in the derivation of \eqref{eq:Nesterov-y} can be relaxed by choosing $a$ in \eqref{eq:V-quad-pd} to be time-varying with the additional condition that $\dot a(t) < 0$.  See also \cite{ross-accelerated-arxiv} for further details on an optimal control theory for accelerated optimization. More specifically, \cite{ross-accelerated-arxiv} shows how new ODEs for new accelerated optimization algorithms can be derived using ideas that generalize the ones introduced in this paper.
\end{remark}

\section*{References}


\begin{thebibliography}{10}

\bibitem{nesterov83}
Yu. E. Nesterov, A method of solving a convex programming problem with convergence rate $\mathcal{O}(1/k^2)$, Soviet Math. Dokl., 27/2 (1983) 371--376 (Translated by A. Rosa).

\bibitem{shi:hi-res-Nes}
B. Shi, S. S. Du, M. I. Jordan, and W. J. Su, Understanding the acceleration phenomenon via high-resolution differential equations, Math. Prog. (2021) https://doi.org/10.1007/s10107-021-01681-8.

\bibitem{alvarez2002}
F. Alvarez, H. Attouch, J. Bolte, P. Redont, A second-order gradient-like dissipative dynamical system with Hessian-driven damping. Applications to optimization and mechanics. J. Math. Pures Appl. 81 (2002) 747--779.

\bibitem{su}
W. Su, S. Boyd, E. J. Candes, A differential equation for modeling Nesterov's accelerated gradient method: theory and insights, J. machine learning research, 17 (2016) 1--43.


\bibitem{boggs71}
P. T. Boggs, The solution of nonlinear system of equations by $A$-stable integration techniques, SIAM J. Numer. Anal. 8/4 (1971) 767--785.

\bibitem{smale76}
S. Smale, A convergent process of price adjustment and global Newton methods, J. mathematical economics, 3 (1976), 107--120.


\bibitem{brown+biggs}
A. A. Brown, M. C. Bartholomew-Biggs, Some effective methods for unconstrained optimization based on the solution of systems of ordinary differential equations, J. optimization theory and applications, 62/2 (1989) 211--224.

\bibitem{gavurin}
M. K. Gavurin, Nonlinear functional equations and continuous analogues of iteration methods, Izv. Vyssh. Uchebn. Zaved. Mat., 5 (1958) 18--31.


\bibitem{rossJCAM-1}
I. M. Ross, An optimal control theory for nonlinear optimization, J. Comp. and Appl. Math., 354 (2019) 39--51.

\bibitem{ross-accelerated-arxiv}
I. M. Ross, An optimal control theory for accelerated optimization, doi = 10.48550/arxiv. 1902.09004, https://arxiv.org/abs/1902.09004,.

\bibitem{ross-book}
I. M. Ross, A Primer on Pontryagin's Principle in Optimal Control, second ed., Collegiate Publishers, San Francisco, CA, 2015.


\bibitem{sontag-book}
E. D. Sontag, Mathematical Control Theory: Deterministic Finite Dimensional Systems, second ed., Springer, New York, NY, 1998.

\bibitem{clarkeLyap}
F. Clarke, Lyapunov functions and feedback in nonlinear control. In: M.S. de Queiroz, M. Malisoff, P. Wolenski (eds) Optimal control, stabilization and nonsmooth analysis. Lecture Notes in Control and Information Science, vol 301. Springer, Berlin, Heidelberg (2004), 267--282.


\bibitem{polyak64}
B. T. Polyak, Some methods of speeding up the convergence of iteration methods, USSR Computational Math. and Math. Phys., 4/5 (1964) 1--17 (Translated by H. F. Cleaves).

\bibitem{polyak-ode}
B. Polyak, P. Shcherbakov, Lyapunov functions: an optimization theory perspective, IFAC PapersOnLine, 50-1 (2017) 7456--7461.



\end{thebibliography}

\end{document}